\documentclass[english,a4paper,12pt,dvipdfm]{amsart}
\usepackage[all,web]{xy}
\usepackage{graphicx}
\usepackage{amsmath} 
\usepackage{amsfonts}
\usepackage{amssymb}
\usepackage{amsthm}
\usepackage{setspace}

\usepackage[linesnumbered, ruled]{algorithm2e}

 \theoremstyle{plain}    
 \newtheorem{thm}{Theorem}[section]
 \numberwithin{equation}{section} 
 \numberwithin{figure}{section} 
 \theoremstyle{remark}   
  
 \theoremstyle{plain}    
 \theoremstyle{plain}    
  
 \theoremstyle{remark} 
 
 \theoremstyle{remark}
 
 \theoremstyle{definition}
  
 \theoremstyle{plain}
 \theoremstyle{plain}    
 \newtheorem{cor}[thm]{Corollary} 
 \theoremstyle{plain}    
 \theoremstyle{definition}
 \newtheorem{defn}[thm]{Definition}
 \theoremstyle{definition}
  \newtheorem{example}[thm]{Example}
 \theoremstyle{plain}    
  
 \theoremstyle{plain}    
 \newtheorem{lem}[thm]{Lemma} 
 \theoremstyle{remark}    
  
 \theoremstyle{remark}    
  
 \theoremstyle{definition}
  
 \theoremstyle{plain}
 \theoremstyle{remark}
 \newtheorem{rem}[thm]{Remark}
 \theoremstyle{remark}

 1

\def\CC{\mathbb{C}}
\def\DD{\mathbb{D}}

\def\QQ{\mathbb{Q}}
\def\RR{\mathbb{R}}

\def\ZZ{\mathbb{Z}}
\def\hyper{\mathcal H}

\newcommand{\mtwo}[4]{\left(
        \begin{matrix}#1&#2\\#3&#4
        \end{matrix}\right)}
\newcommand{\abcd}[4]{\left(
        \begin{smallmatrix}#1&#2\\#3&#4\end{smallmatrix}\right)}

\def\inner<#1,#2>{{\left\langle{{#1},{#2}}\right\rangle}}
\def\sl2of#1{\textrm{SL}_2(#1)}

\def\inner<#1,#2>{{\left\langle{{#1},{#2}}\right\rangle}}

\usepackage{babel}

\begin{document}

\title{A level $N$ reduction theory of indefinite binary quadratic forms}
\author{Carlos Casta\~no-Bernard}
\email{ccastanobernard@gmail.com}
\urladdr{http://sarva.no-ip.info/ccb/index.html}

\begin{abstract}
In this paper we study a geometric coding algorithm for
indefinite binary quadratic forms $Q$ for
the congruence subgroup $\Gamma^0(N)$,
with respect to the usual fundamental domain $\mathcal{F}_N$,
where $N$ is assumed prime.
The cycles $Q_1,\dots,Q_n$ that this algorithm produces are such that
the the corresponding paths $\gamma_1,\dots,\gamma_n$
in the Riemann surface $X_0(N)(\CC)$ have a nice behaviour around
the elliptic points of order $2$.
\end{abstract}

\maketitle

\pagenumbering{roman}
\setcounter{page}{0}



\tableofcontents


\pagenumbering{arabic}
\pagestyle{headings}

%
%
%
\section{Introduction}
In this paper we study a geometric coding algorithm for
indefinite binary quadratic forms $Q$ for
the congruence subgroup $\Gamma^0(N)$,
with respect to the usual fundamental domain $\mathcal{F}_N$,
where $N$ is assumed prime.
(Cf. Schoeneberg book~\cite{schoeneberg:ellipticmodular}.)
Following Ogg~\cite{ogg:real}
let us call $E$-points the elliptic points of order $2$
with respect to the action of $\Gamma^0(N)$
on the upper half-plane
\begin{displaymath}
\hyper=\{\tau\in\CC\colon\Im(\tau)>0\}.
\end{displaymath}
The cycle
\begin{displaymath}
Q_1,\dots,Q_n
\end{displaymath}
that this algorithm attaches to $Q$ deals with the $E$-points
is such that the the corresponding paths
\begin{displaymath}
\gamma_1,\dots,\gamma_n
\end{displaymath}
define in a natural way
a regular path on the  Riemann surface $X_0(N)(\CC)$.

The motivation to write this paper was mostly utilitarian.
In a forthcoming paper of the author
this algorithm is applied to the computation of the
homology class of each of the connected components of
the real locus of the Atkin-Lehner quotient $X_0^+(N)$
of $X_0(N)$ associated to the Fricke involution $w_N$.

\section{The level $1$ case}

\subsection{Preliminaries}
Let $[A,B,C](X,Y)=AX^2 + BXY + CY^2$ denote a binary
quadratic forms with real coefficients $A$, $B$, and $C$.
Note that for any real numbers
$\alpha, \beta, \gamma$ and 
$\delta$ we have
\begin{equation*}
[A,B,C](\alpha X + \beta Y, \gamma X + \delta Y)=
\left(
\begin{matrix}
\alpha^2&\alpha\gamma&\gamma^2\\
2\alpha\beta&\alpha\delta+\beta\gamma&2\gamma\delta\\
\beta^2&\beta\delta&\delta^2\\
\end{matrix}
\right)
\left(
\begin{matrix}
A\\
B\\
C\\
\end{matrix}
\right).
\end{equation*}
Let 
$M=\abcd{\alpha}{\beta}{\gamma}{\delta}\in\textrm{GL}_2(\RR)$ and define
\begin{equation*}
\sigma_M=
\det (M)
\left(
\begin{matrix}
\alpha^2&\alpha\gamma&\gamma^2\\
2\alpha\beta&\alpha\delta+\beta\gamma&2\gamma\delta\\
\beta^2&\beta\delta&\delta^2\\
\end{matrix}
\right)^{-1}.
\end{equation*}
Let $\mathcal{Q}$ denote the set of binary quadratic forms.
Note that
\begin{equation*}
\langle Q, Q^\prime\rangle=BB^\prime - 2(AC^\prime + A^\prime C).
\end{equation*}
makes $\mathcal{Q}$ a quadratic space and
$M\mapsto \sigma_M$ defines a group homomorphism
\begin{equation*}
\sigma:\textrm{GL}_2(\RR)\longrightarrow O^+(\mathcal{Q})
\end{equation*}
of the general linear group GL$_2(\RR)$ into
the group of orientation preserving automorphisms
$O^+(\mathcal{Q})$ of 
the ternary quadratic form $|Q|^2=\langle Q,Q\rangle$
on $\mathcal{Q}$.\footnote{
In fact the homomorphism $\sigma_M$ restricted to
SL$_2^\pm(\RR)=\{m\in\textrm{GL}_2(\RR)\colon\det m=\pm 1\}$
is a degree $2$ cover of SL$_2(\RR)$
known as the spin representation attached to the
ternary quadratic form $|Q|^2=\langle Q,Q\rangle$.}
For each matrix $M\in\textrm{GL}_2(\RR)$ and each binary quadratic form
$Q\in\mathcal{Q}$ we write $M\circ Q=\sigma_MQ$.
The ternary quadratic form $|Q|^2=\langle Q,Q\rangle$ is called the
\textit{discriminant} of the binary quadratic form $Q=[A,B,C]$.
We say $Q$ is
\textit{definite} (resp. \textit{indefinite})
if $|Q|^2>0$ (resp. $|Q|^2<0$).
If $Q$ is \textit{definite} and 
$A>0$ then we say $Q$ is \textit{positive definite}.

\begin{defn}\label{defn:circ}
Two binary quadratic forms 
$Q_1$ and $Q_2$ are called \textit{equivalent} if there exists 
$M\in\textrm{SL}_2(\ZZ)$ such that $M\circ Q_1 = Q_2$.
\end{defn}

\begin{example}
Any form $Q = [A, B, C]$ is equivalent to 
\begin{equation*}
S\circ Q = [C, -B, A],
\end{equation*}
where $S = \abcd{0}{-1}{1}{0}\in\textrm{SL}_2(\ZZ)$.
\end{example}

\begin{example}
Any form $Q = [A, B, C]$ is equivalent to 
\begin{equation*}
T^t\circ Q = [A, -2tA + B, t^2A -tB + C],
\end{equation*}
where $T = \abcd{1}{t}{1}{0}\in\textrm{SL}_2(\ZZ)$,
for each $t\in\ZZ$.
\end{example}

\begin{defn}
We say a positive definite binary quadratic form
$Q=[A,B,C]$ is \textit{normal} if $-A \leq -B < A$.
\end{defn}

Every positive definite form $Q$ is equivalent to a normal form, namely
$Q^\textit{nrm}=T^t\circ Q$, where
\begin{equation*}
t=-\left\lfloor \frac{A - B}{2A} \right\rfloor
\end{equation*}
and $\lfloor x \rfloor$ denotes the floor of a real number $x$.
The form $Q^\textit{nrm}$ is known as the \textit{normalisation} of $Q$.

\begin{defn}
We say a positive definite binary quadratic form $Q=[A,B,C]$
is \textit{reduced} if $Q$ is normal and if
any of the following two conditions is true.
\begin{enumerate}
\item $A < C$
\item $A = C$ and $-B \leq 0$ 
\end{enumerate}
\end{defn}

\begin{algorithm}\label{alg:gauss}
  \KwData{any positive definite binary quadratic form $Q$}
  \KwResult{a reduced form $Q^\textit{red}$ equivalent to $Q$}
  $Q\longleftarrow Q^\textit{nrm}$\;
  \While{not reduced}{
    $Q\longleftarrow (S\circ Q)^\textit{nrm}$\;
    }
  \caption{The reduction of a positive definite form $Q$}
\end{algorithm}

\begin{figure}
\begin{center}
\includegraphics{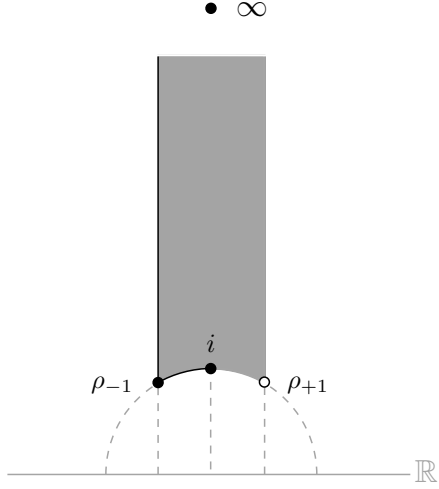}
\end{center}
\caption{The fundamental set $\mathcal{F}^\prime$.}
\label{cap:fundom}
\end{figure}

The definition of reduced form
may be rephrased in geometric terms as follows.
To each positive definite form $Q$ attach
the point $\pi_Q\in\hyper$ defined by
\begin{equation*}
\pi_Q=\frac{-B + |Q|}{2A},
\end{equation*}
where $|Q|$ denotes the square root of 
$|Q|^2$ such $\Im(|Q|)>0$.
So $Q$ is reduced if and only if
$\pi_Q\in\mathcal{F}^\prime$,
where 
\begin{equation*}
\mathcal{F}^\prime=
\left\{
\tau\in\mathcal{H}\colon
-\frac{1}{2} \leq \Re(\tau) <\frac{1}{2}
\right\}
\cap
\left\{
\tau\in\mathcal{H}\colon \, |\tau|>1\,\textit{or}\,\,
(|\tau|=1\,\,\textit{and}\,\,\Re(\tau)\leq 0)
\right\},
\end{equation*}
which is depicted by Figure~\ref{cap:fundom}.
If we define $\tau\mapsto M\tau$ for each
$M=\abcd{\alpha}{\beta}{\gamma}{\delta}\in\textrm{SL}_2(\RR)$
by the M\"obius transformation
\begin{equation*}
M\tau=\frac{\alpha\tau +\beta}{\gamma\tau + \delta},
\end{equation*}
then we have an action of $\textrm{SL}_2(\RR)$ on $\hyper$ such that
$M \pi_Q=\pi_{M\circ Q}$.

While performing Algorithm~\ref{alg:gauss} on a given 
positive definite form $Q$ we may obtain a product of matrices
\begin{equation*}
R_Q=T^{n_l} S \cdots S T^{n_2} S T^{n_1}\in \textrm{SL}_2(\ZZ)
\end{equation*}
such that $R_Q\circ Q$ is reduced,
where $S=\abcd{0}{-1}{1}{0}$ and $T=\abcd{1}{1}{0}{1}$.
Note that in terms of the above geometric setting
we have $R_Q\pi_Q\in\mathcal{F}$.
So given any $\tau\in\hyper$ there is a matrix $R_\tau\in\textrm{SL}_2(\ZZ)$
such that $R_\tau \tau\in\mathcal{F}$.
In fact $R_\tau$ is uniquely determined by $\tau$,
except when $\tau$ is equivalent to either
$i=\sqrt{-1}$, with stabiliser of order $2$
or $\rho_{-1}=\frac{-1+\sqrt{-3}}{2}$, with stabiliser of order $3$.

\begin{rem}
Later,
in connection with the coding algorithm
we will find convenient to use
the closure $\mathcal{F}_N=(\mathcal{F}^\prime)^\textit{clo}$ of
the fundamental set $\mathcal{F}_N^\prime$.
In the literature the closed set $\mathcal{F}\subset\hyper$
is known as a \textit{fundamental region}.
\end{rem}

%
\subsection{Indefinite forms}\label{sub:redindef}
Let $\textit{sign}(x)$ denote the sign of a real number $x$.
We will find convenient to assume
the somewhat non-standard convention that $\textit{sign}(0) = 1$,
as opposed to the usual $\textit{sign}(0) = 0$.

\begin{defn}\label{defn:redindef}
We say an indefinite binary quadratic form
$Q=[A,B,C]$ is \textit{reduced}\footnote{
A closely related definition of reduced indefinite form
may be found in Choie and Parson's article~\cite{choie:rational},
where $Q=[A,B,C]$ is called reduced
if $A > 0$, $C > 0$, and $A + C < B$.}
if $B>0$ and if any of the following two conditions is satisfied.
\begin{description}
\item[R1] $2|A+C| < B$
\item[R2] $2(A+C) = -\textit{sign}(A) B$
\end{description}
\end{defn}

From now on let us assume $Q=[A,B,C]$ is an indefinite binary quadratic form
with $A,B$ and $C$ in $\ZZ$.
We shall give a geometric interpretation of
the concept of reduced indefinite form.
Provide the upper half-plane $\hyper$ with the hyperbolic metric
\begin{displaymath}
ds=\frac{|d\tau|}{y},
\end{displaymath}
where $\tau=x+iy$,
with $x$ and $y$ the usual (real-valued) coordinate functions of $\CC$.
Let $\{\tau_1,\tau_2\}$ denote the (oriented) geodesic that joins a
given point $\tau_1\in\hyper^*$ to a given point $\tau_2\in\hyper^*$.
The \textit{geodesic} $\gamma_Q$ attached to an indefinite form $Q$ is
\begin{displaymath}
\gamma_Q=\left\{
\begin{array}{ll}
\{\frac{-B -|Q|}{2A},\frac{-B +|Q|}{2A}\}&\textrm{if $A\not=0$}\\ 
\{i\infty,-\frac{C}{B}\}&\textrm{if $A=0$ and $B>0$}\\
\{-\frac{C}{B},i\infty\}&\textrm{if $A=0$ and $B<0$}\\
\end{array}
\right.
\end{displaymath}
where $|Q|$ denotes the positive square root of $|Q|^2$.
We also write $\sqrt{D}=|Q|$, where $D=B^2 - 4AC$.
The \textit{intersection number} 
$I(\gamma_1,\gamma_2)$ of
given (oriented) geodesics arcs $\gamma_1$ and $\gamma_2$
on the upper half-plane $\hyper^*$ is defined by
\begin{displaymath}
I(\gamma_1,\gamma_2)=\left\{
\begin{array}{rl}
 1,& \textrm{if $\,\#(\gamma_1\cap\gamma_2) = 1 $ and
$\{\gamma^\prime_1(P),\gamma^\prime_2(P)\}$ has {\it positive} orientation}\\
-1,& \textrm{if $\,\#(\gamma_1\cap\gamma_2) = 1 $ and
$\{\gamma^\prime_1(P),\gamma^\prime_2(P)\}$ has {\it negative} orientation}\\
0, &\textrm{otherwise.}\\
\end{array}
\right.
\end{displaymath}
where $\gamma^\prime(P)$ denotes the unit tangent vector of a geodesic
$\gamma$ at a point $P\in\gamma$. 

\begin{example}
Suppose $Q=[0,1,0]$ and $U=[1,0,-1]$.
In fact $\gamma_Q=\{i\infty,0\}$ and $\gamma_U=\{-1,1\}$.
Clearly $\gamma_Q$ and $\gamma_U$ meet at $P=\pi_{[1,0,1]}=i$
(and no other point).
In particular $I(\gamma_Q,\gamma_U)=\pm 1$.
To determine the sign we note that the
tangent vector of $\gamma_Q$ at $P$ is given by
$\gamma^\prime_Q(P)=(0,-1)$,
and the tangent vector of $\gamma_U$ at $P$ is given by
$\gamma^\prime_U(P)=(1,0)$.
Since the determinant
\begin{equation*}
\left|
\begin{matrix}
0&-1\\
1&0\\
\end{matrix}
\right|>0
\end{equation*}
we have
$I(\gamma_Q,\gamma_U)=1$.
\end{example}

Put $\sigma=\{\rho_{-1},\rho_{+1}\}$,
where as before $\rho_{\pm 1} = \frac{\pm 1 +\sqrt{-3}}{2}$.
Note $\sigma\subset\gamma_U$, where $U=[1,0,-1]$.
To simplify our exposition assume from now on that
\begin{equation*}
\fbox{$D$ \it is not a perfect square.}
\end{equation*}

\begin{lem}
An indefinite form $Q$ is reduced if and only if
$I(\gamma_Q,\sigma)=1$
and
$\gamma_Q\cap\mathcal{F}^\textit{int}\not=\emptyset$.
\end{lem}

\begin{proof}
First note that given two binary quadratic forms $Q_1$ and $Q_2$,
where $Q_1$ is positive definite and $Q_2$ is indefinite,
the condition that the point $\pi_{Q_1}$ lies on
the geodesic $\gamma_{Q_2}$
is equivalent to $\langle Q_1, Q_2\rangle=0$.
So $\gamma_Q\cap\gamma_U\not=\emptyset$ means there is 
a positive definite form $Q^\prime$ such that
\begin{enumerate}
\item $\langle Q^\prime,U\rangle=0$
\item $\langle Q^\prime,Q\rangle=0$
\end{enumerate}
Write $Q^\prime=[A^\prime,B^\prime,C^\prime]$.
Clearly (1) and (2) hold if and only if
\begin{equation*}
A^\prime=C^\prime \quad\textrm{and}\quad B^\prime B-2(AC^\prime +A^\prime C)=0.
\end{equation*}
This means that $B^\prime B-2(AA^\prime +A^\prime C)=0$.
In other words
\begin{equation*}
\frac{B^\prime}{2A^\prime} = \frac{A+C}{B}.
\end{equation*}
Thus $\gamma_Q$ meets $\gamma_U$ at a point $\pi_{Q^\prime}$
on the geodesic arc $\sigma=\{\rho_{-1},\rho_{+1}\}$ if and only if
\begin{equation*}
2|A+C|\leq |B|.
\end{equation*}
Now we split the proof in two cases
according as $\pi_{Q^\prime}$ is an endpoint of $\sigma$ or not.

\noindent\textbf{Case I:} Suppose
\begin{equation*}
2(A+C)=\mp B.
\end{equation*}
This condition is equivalent to $\rho_{\pm 1}\in\gamma_Q$.
So $\gamma_{Q}\cap\mathcal{F}^\textit{int}\not=\emptyset$ if and only if
\begin{equation*}
\mp\frac{B}{2A}<0,
\end{equation*}
and $I(\gamma_Q,\sigma)=1$ if and only if $\pm A>0$.
In other words, $2(A + C)=-\textit{sign}(A)B$ and $B>0$.

\noindent\textbf{Case II:} Now suppose
\begin{equation*}
2|A+C| < |B|.
\end{equation*}
Since the point $\pi_{Q^\prime}$ is not an endpoint of $\sigma$,
we have $\gamma_Q\cap\mathcal{F}^\textit{int}\not=\emptyset$.
Moreover, the condition $B>0$ is equivalent to $I(\gamma_Q,\sigma)=1$.
The lemma follows.
\end{proof}

We will shortly introduce Algorithm~\ref{alg:reduce},\footnote{
The author owes the main ideas underlying this algorithm to
a talk given by Richard Pinch.}
which is a reduction algorithm of
indefinite binary quadratic forms that we will use
in connection with Algorithm~\ref{alg:cycle}.
The latter algorithm is
a variant on the geometric coding algorithm of Katok~\cite{katok:coding}.

We call an indefinite form $Q=[A,B,C]$ \textit{nearly reduced}
if  $3A^2 \leq D$, where $D=B^2 - 4AC$ is the discriminant of $Q$.
To each nearly reduced indefinite form $Q$ we attach the interval
\begin{displaymath}
J_Q=\left\{
\begin{array}{ll}
\left[t_{-1}^\prime, t_{+1}^\prime\right],&
\textit{if $D > 4A^2$}\\
&\\
\left[t_{-\textit{sign}(A)}, t_{-\textit{sign}(A)}^\prime\right],&
\textit{if $D < 4A^2$}\\
\end{array}
\right.
\end{displaymath}
where
\begin{equation*}
t_{\pm 1}=\pm\frac{1}{2} + \frac{B + \sqrt{D-3A^2}}{2A},
\quad\textrm{and}\quad
t_{\pm 1}^\prime=\pm\frac{1}{2} +\frac{B - \sqrt{D-3A^2}}{2A}.
\end{equation*}
Note that the assumption that $Q$
is nearly reduced guarantees that the end points of the interval $J_Q$ are real.
 
\begin{lem}\label{lem:reduce}
Suppose $Q$ is a nearly reduced indefinite form.
We have $I(\gamma_{T^t\circ Q},\sigma)=1$
if and only if $t\in J_Q$,
where $T^t=\abcd{1}{t}{0}{1}$.
\end{lem}

\begin{proof}
First we define an auxiliary function $f_Q$ attached to
the nearly reduced form $Q=[A,B,C]$ by letting
\begin{equation*}
f_Q(t)=\left|\frac{A_t+C_t}{B_t}\right|,
\end{equation*}
where
\begin{displaymath}
\begin{array}{lll}
A_t&=&A,\\
B_t&=&B-2At,\\
C_t&=&At^2 -Bt +C,\\
\end{array}
\end{displaymath}
with $t\in\RR$.
Note that $\gamma_{T^t\circ Q}$
has non-trivial intersection with the geodesic arc
$\sigma=\{\rho_{-1},\rho_{+1}\}$ if and only if $t\in\RR$ is such that
\begin{equation*}
f_Q(t)\leq \frac{1}{2}.
\end{equation*}
We will see that the set of such $t$ may be expressed as
a union of two closed intervals.
Our first step will be to find the endpoints of these intervals,
i.e. the solutions $t\in\RR$ of $f_Q(t)=\frac{1}{2}$.
The latter is equivalent to any of the following equations
\begin{equation*}
\begin{split}
2|A_t+C_t|&=|B_t|,\\
2|At^2 - Bt + A + C|&=|B - 2At|,\\
2(At^2 - Bt + A + C)&=\pm (B - 2At).
\end{split}
\end{equation*}
A straight forward calculation shows that the four real numbers
$t_{+1}$, $t_{-1}$, $t_{+1}^\prime$, and $t_{-1}^\prime$,
defined above are all the solutions of $f_Q(t)=\frac{1}{2}$.
Now the proof splits in two cases.

\noindent\textbf{Case I:} Suppose $D > 4A^2$.
This condition means that the intervals
$[t_-, t_+]$ and $[t_-^\prime, t_+^\prime]$ are disjoint.
It is straight forward to see that $f_Q(t)\leq \frac{1}{2}$,
if and only if $t\in [t_-, t_+]\cup [t_-^\prime, t_+^\prime]$,
with $f_Q(t) = \frac{1}{2}$ exactly at the endpoints.
Clearly $I(\gamma_{T^t\circ Q},\sigma)=1$ if and only if
\begin{equation*}
t\in\left[t_{-1}^\prime,t_{+1}^\prime\right].
\end{equation*}
\noindent\textbf{Case II:} Suppose $D < 4A^2$.
This condition means that the intervals
$[t_-, t_-^\prime]$ and $[t_+, t_+^\prime]$ are disjoint.
Again it is straight forward to see that $f_Q(t)\leq \frac{1}{2}$,
if and only if $t\in [t_-, t_-^\prime]\cup[t_+, t_+^\prime]$,
with $f_Q(t) = \frac{1}{2}$ exactly at the endpoints.
Note that $I(\gamma_{T^t\circ Q},\sigma)=1$ if and only if
\begin{equation}\label{eqn:gammarho}
\left|\rho_{-\textit{sign}(A)} -t + \frac{B}{2A}\right|\leq \frac{\sqrt{D}}{2|A|}.
\end{equation}
By squaring the left-hand side of Inequality~\ref{eqn:gammarho} we get
\begin{equation*}
\left|\rho_{-\textit{sign}(A)} -t + \frac{B}{2A}\right|^2 =
\left(\rho_{-\textit{sign}(A)} -t + \frac{B}{2A}\right)
\left(\bar{\rho}_{-\textit{sign}(A)} -t + \frac{B}{2A}\right)=
\end{equation*}
\begin{equation*}
\rho_{-\textit{sign}(A)}\bar{\rho}_{-\textit{sign}(A)} +
\left(\rho_{-\textit{sign}(A)} +
\bar{\rho}_{-\textit{sign}(A)}\right)\left(-t+\frac{B}{2A}\right) +
t^2 - \frac{B}{A}t + \frac{B^2}{4A^2}=
\end{equation*}
\begin{equation*}
1 -\textit{sign}(A)\left(-t + \frac{B}{2A}\right) +
t^2 - \frac{B}{A}t + \frac{B^2}{4A^2}
\end{equation*}
So Inequality~\ref{eqn:gammarho} is equivalent to
\begin{equation*}
t^2 +\left(\textit{sign}(A)-\frac{B}{A}\right)t -\textit{sign}(A)\frac{B}{2A} +
\frac{B^2}{4A^2} -\frac{D}{4A^2} + 1\leq 0.
\end{equation*}
After a straight forward calculation we may see that
$t\in\RR$ satisfies Inequality~\ref{eqn:gammarho} if and only if
$t\in \left[t_{-\textit{sign}(A)},t_{-\textit{sign}(A)}^\prime\right]$,
and the lemma follows.
\end{proof}

\begin{cor}\label{cor:reduce}
If $Q=[A,B,C]$ is an indefinite form such that
$\gamma_Q\cap \mathcal{F}\not=\emptyset$,
then there is exactly one $t\in\ZZ$ such that $T^t\circ Q$ is reduced,
unless $D < 4A^2$ and either $\rho_{-1}$ or $\rho_{+1}$ lies in $\gamma_Q$.
\end{cor}

\begin{proof}
We split the proof splits in two cases as follows.

\noindent\textbf{Case I}
Suppose $Q=[A,B,C]$ is an indefinite form such that $D > 4A^2$.
Note that the closed interval $J_Q$ has length $1$.
Now we consider two sub-cases.
\begin{description}
\item[Sub-case I] There is exactly one point $t\in J_Q\cap \ZZ$.
Note that
$\gamma_{T^t\circ Q}\cap \mathcal{F}^\textit{int}\not=\emptyset$.
Thus by Lemma~\ref{lem:reduce} we may see $T^t\circ Q$ is reduced.
\item[Sub-case II] The endpoints $t_{-1}^\prime$ and $t_{+1}^\prime$ are integral.
Clearly $t=t_{\textit{sign}(A)}^\prime$ is the only point such that 
$\gamma_{T^t\circ Q}\cap\mathcal{F}^\textit{int}\not=\emptyset$,
and again by applying Lemma~\ref{lem:reduce}
we may see that $T^t\circ Q$ is reduced.
\end{description}

\noindent\textbf{Case II}
Suppose $Q=[A,B,C]$ is such that $D < 4A^2$,
but neither $\rho_{-1}$ nor $\rho_{+1}$ lies in $\gamma_Q$.
Clearly $\gamma_{T^t\circ Q}\cap\mathcal{F}^\textit{int}\not=\emptyset$.
Now using the Lemma~\ref{lem:reduce} we may see that
$I(\gamma_{T^t\circ Q},\sigma)=1$ for suitable $t\in\{-1,0,1\}$.
The corollary follows.
\end{proof}

\begin{lem}\label{lem:rhoint}
Suppose $Q$ is an indefinite form such that
either $\rho_{-1}\in \gamma_{Q}$ or $\rho_{+1}\in \gamma_{Q}$.
Then 
$\gamma_{Q}\cap\mathcal{F}^\textit{int}=\emptyset$ if and only if
$\gamma_{S\circ Q}\cap\mathcal{F}^\textit{int}\not=\emptyset$.
\end{lem}

\begin{proof}
On the one hand note that
the assumption $\rho_{\pm 1}\in\gamma_{Q}$ is equivalent to
\begin{equation}\label{eqn:rho01}
2(A + C)=\mp B,
\end{equation}
where $Q = [A, B, C]$.
On the other hand note 
$\gamma_{Q}\cap \mathcal{F}^\textit{int}=\emptyset$
is equivalent to
\begin{equation}\label{eqn:rho02}
-\frac{B}{2A}\in (0,\pm 1).
\end{equation}
Using Equation~\ref{eqn:rho01} we may see
Inclusion~\ref{eqn:rho02} is equivalent to
\begin{equation*}
\frac{C}{A}\in (0, - 1).
\end{equation*}
An immediate consequence of the latter is that the form
\begin{equation*}
Q^\prime=S\circ Q=[C, -B, A]
\end{equation*}
satisfies
$\gamma_{Q^\prime}\cap \mathcal{F}^\textit{int}\not=\emptyset$.
The converse is clearly true and the lemma follows.
\end{proof}

\begin{defn}\label{defn:nrm}
We say an indefinite form $Q$ is \textit{normalisable}
if there is $\delta(Q)\in\ZZ$ such that $T^{\delta(Q)}\circ Q$ is reduced.
Such form $Q^\textit{nrm} = T^{\delta(Q)}\circ Q$ is called its
\textit{normalisation}.
\end{defn}

\begin{defn}
The \textit{tip} $\widehat{\gamma}$ of a geodesic $\gamma\subset\hyper$
is the point
\begin{equation*}
\widehat{\gamma}=\frac{x+y}{2}+\left|\frac{x-y}{2}\right|i\in\gamma
\end{equation*}
where $\gamma=\{x,y\}$ and $x$ and $y$ in $\RR$.
We define the \textit{tip} $\widehat{Q}$ of an indefinite form $Q$ as
the positive definite binary quadratic form $P$ such that
$\pi_{P}=\widehat{\gamma_Q}$.
If $Q$ is integral we assume $\widehat{Q}$ is primitive.
\end{defn}

We will find it convenient
to consider pointed spaces $(\gamma_Q,\tau_0)$,
with $\tau_0\in\gamma_Q$,
e.g. $\tau_0=\widehat{\gamma_Q}$.
By a slight abuse of notation
we sometimes write $(Q,\tau_0)$ instead of $(\gamma_Q,\tau_0)$.

\begin{algorithm}\label{alg:reduce}
  \KwData{$(Q,\tau_0)$ with $\tau_0\gamma_Q$}
  \KwResult{reduced form $Q^\textit{red}$ equivalent to $Q$}
$Q\longleftarrow R_{\tau_0} \circ Q$\;
\eIf{$\gamma_Q\cap\mathcal{F}^\textit{int}\not=\emptyset$
}{$Q\longleftarrow Q^{\textit{nrm}}$
}{$Q\longleftarrow (S\circ Q)^{\textit{nrm}}$
}
  \caption{The reduction of $Q$}
\end{algorithm}

Now we prove that Algorithm~\ref{alg:reduce} is correct.
First note that
the correctness of Algorithm~\ref{alg:gauss} (due to Gau\ss )
implies that after Step 1
the geodesic $\gamma_Q$ and 
the (closed) fundamental region $\mathcal{F}$ share at least one point.
If also the interior $\mathcal{F}^\textit{int}$ and $\gamma_Q$
have non-trivial intersection,
then Corollary~\ref{cor:reduce} implies that
the normalisation $Q^\textit{nrm}$ of $Q$ is reduced.
However,
it may sometimes happen that
$\mathcal{F}^\textit{int}\cap\gamma_Q=\emptyset$,
which means that
$\gamma_Q$ contains one of the lower vertices of
the fundamental region $\mathcal{F}$.
But Lemma~\ref{lem:rhoint} tells us that
the transformed form $S\circ Q$ now belongs to the above case,
and thus $(S\circ Q)^\textit{nrm}$ is reduced.
Therefore Algorithm~\ref{alg:reduce} is correct.

%
\subsection{Regular path attached to an indefinite form}\label{sub:reg}
As before let $Q=[A, B, C]$ be primitive indefinite form
with non-square discriminant $D = B^2 - 4AC$ and define
\begin{equation*}
M_Q = \mtwo{x - By}{-2Cy}{2Ay}{x + By}\in\Gamma(1),
\end{equation*}
where $(x, y)\in\ZZ\times\ZZ$ is the
fundamental solution to Pell's equation
\begin{equation}\label{eqn:pell}
X^2 - DY^2 = 1
\end{equation}
such that both, $x > 0$ and $y > 0$,
so that the eigenvalue
$\lambda_Q=x+y\sqrt{D}\in\mathcal{O}_D^\times$
of $M_Q$ satisfies $\lambda_Q>1$.
Note that $M_Q$ is an automorphism of the quadratic form $Q$.
So given a base-point $\tau_0\in\gamma_Q$,
the geodesic segment
$\{\tau_0,M_Q\tau_0\}$ is a subset of the geodesic
\begin{displaymath}
\gamma_Q=\left\{\frac{-B-\sqrt{D}}{2A},\frac{-B+\sqrt{D}}{2A}\right\}.
\end{displaymath}
Now,
by standard results of hyperbolic geometry
as in Beardon's book~\cite{beardon:discrete},
we may see that the arc-length parametrisation of
the (oriented) geodesic segment $\{\tau_0, M_Q\tau_0\}$
is the restriction of
the arc-length parametrisation of
the (oriented) geodesic $\gamma_Q$ with base point $\tau_0$
to the closed interval $[0,l(Q)]\subset\RR$,
where $l(Q)=2\log(\lambda_Q)$ is
the hyperbolic length of $\{\tau_0,M_Q\tau_0\}$.
Now let
\begin{equation*}
g_{Q,\tau_0}\colon \left[0, l(Q)\right]\longrightarrow X(1).
\end{equation*}
be the parametrisation induced in the quotient $X(1)$
by the above restriction.
Note that
$(\gamma_Q)_{\Gamma(1)}=\{\tau_0,M_Q\tau_0\}_{\Gamma(1)}$,
where $X_{\Gamma(1)}$ denotes the image of
a subset $X\subset\hyper^*$ in
the quotient $X(1)=\Gamma(1)\backslash\hyper^*$.
Since $Q$ is primitive,
the assumption that $(x,y)$ is a fundamental solution to
Equation~\ref{eqn:pell} is equivalent to the assertion that 
$M_Q$ is not a non-trivial power of an element in
the modular group $\Gamma(1)$.
So \textit{essentially} $g_{Q,\tau_0}$ traces out its image
$(\gamma_Q)_{\Gamma(1)}=\{\tau_0,M_Q\tau_0\}_{\Gamma(1)}$
only once.
We call $g_{Q,\tau}$
the \textit{closed geodesic} associated to
the form $Q$ and the point $\tau_0\in\gamma_{Q}$.\footnote{
In the literature the closed geodesic $g_{Q,\tau_0}$
is also known as a \textit{prime geodesic}.
By the Prime Geodesic Theorem of Sarnak~\cite{sarnak:thesis}
we know that the distribution of the lengths $l(Q)$ is similar to
the distribution of prime numbers.}

Let $X$ be any differentiable manifold
and let $T$ be a positive real number.
A differentiable path $g:[0,T]\longrightarrow X$ is called
\textit{regular} if the derivative $g^\prime(t)\not=0$,
for each $t\in[0,T]$.

\begin{lem}\label{lem:reg}
Again let $Q$ be an indefinite form
and $\tau$ a point of $\gamma_Q$.
If the image of the path $\gamma_Q\subset\hyper$
contains an elliptic point,
then the path $g_{Q,\tau}\subset\hyper$ is not regular,
else the path $g_{Q,\tau}\subset\hyper$ is regular.
\end{lem}

\begin{proof}
It is well-known that
the complex-analytic structure of
the affine modular curve $Y(1)=\Gamma(1)\backslash\hyper$ is
induced by the map $j\colon\hyper\mapsto\CC$,
where $j(\tau)$ is the elliptic modular function
\begin{displaymath}
j(\tau)=\frac{1}{q} + 744 + 196884q+\dots,
\end{displaymath}
(See Knapp's book~\cite{knapp:elliptic}.)
Moreover,
given any point $b\in\hyper$ there is
a sufficiently small neighbourhood $U$ (resp. $V$) of $b$ (resp. $j(b)$)
together with
an holomorphic identification $\phi\colon U\longrightarrow \DD$
(resp. $\psi\colon V\longrightarrow \DD$)
with the unit disk
\begin{displaymath}
\DD=\{z\in\CC\colon |z|<1\}
\end{displaymath}
such that the diagram
\begin{equation*}
\xymatrix{
U \ar[d]^{\phi}\ar [r]^j& V\ar[d]^{\psi}\\
\DD\ar [r]& \DD\\
}
\end{equation*}
\begin{equation*}
\xymatrix{
z \ar @{|->}[r]& z^n\\
}
\end{equation*}
commutes, where
\begin{displaymath}
n=\left\{
\begin{array}{ll}
3,&\textit{if $\,b$ is an elliptic point of order $3$,}\\
2,&\textit{if $\,b$ is an elliptic point of order $2$,}\\
1,&\textit{otherwise.}\\
\end{array}
\right.
\end{displaymath}
So there is $t_0\in[0,l(Q)]$ such that $g_{Q,\tau_0}(t_0)$
is an elliptic point then
if and only if the derivative $g_{Q,\tau_0}^\prime(t_0)=0$.
The lemma follows.
\end{proof}

For the rest of the subsection
suppose $g_{Q,\tau_0}$ fails to be regular
at a point $t_0\in[0,T]$,
i.e. $g_{Q,\tau_0}(t_0)=\tau_e$,
where $\tau_e$ is (the image in $X(1)$ of) an elliptic point
of order $e=2$ or $3$.
Shortly we shall see how,
under certain mild conditions,
a suitable restriction of the path $g_{Q,\tau_0}$
may be extended in a natural way
to produce an analytic regular path.
But first we need to prove the following
two elementary lemmas on square roots of
certain elements of the group $\textrm{SL}_2(\RR)$,
and also introduce a basic definition
immediately after the lemmas.
The mentioned regular path is actually constructed
in the proof of Theorem~\ref{thm:cont}.

\begin{lem}\label{lem:sqr}
Let $Q$ and $M_Q$ be as above.
Also let $\epsilon=x^\prime+y^\prime\sqrt{D}\in\mathcal{O}_D^\times$
(with $x^\prime$ and $y^\prime$ in $\frac{1}{2}\ZZ$)
be a fundamental unit of
the real quadratic order $\mathcal{O}_D$ of discriminant $D$.
Assume that $\epsilon$ has negative norm $\mathcal{N}(\epsilon)=-1$
and also that $\epsilon>1$.
Then
\begin{equation*}
M_Q^{\frac{1}{2}}=\frac{1}{\sqrt{D}}
\mtwo{Dy^\prime-Bx^\prime}{-2Cx^\prime}{2Ax^\prime}{Dy^\prime+Bx^\prime}
\in\textrm{SL}_2(K)
\end{equation*}
is a square root of $M_Q$ such that
the arc-length parametrisation of
the geodesic segment $\{\tau_0,M_Q^\frac{1}{2}\tau_0\}$
is the restriction of
the arc-length parametrisation of
the geodesic segment $\{\tau_0,M_Q\tau_0\}$
to the interval $[0,\frac{1}{2}l(Q)]$.
\end{lem}

\begin{proof}
Consider the matrix
\begin{equation*}
M=\mtwo{\frac{-B+\sqrt{D}}{2A}c}{\frac{-B-\sqrt{D}}{2A}d}{c}{d}
\in\textrm{SL}_2(K)
\end{equation*}
where $c$ and $d$ are integers such that $cd=\frac{A}{\sqrt{D}}$.
It is plain that the matrix $M$
maps the elements $0$ and $\infty$ as follows:
\begin{equation*}
\begin{split}
0 \mapsto \frac{-B-\sqrt{D}}{2A},\\
\infty \mapsto \frac{-B+\sqrt{D}}{2A}.\\
\end{split}
\end{equation*}
Now recall that
the solution $(x,y)\in\ZZ\times\ZZ$ of Equation~\ref{eqn:pell}
is minimal subject to the conditions $x>0$ and $y>0$.
In particular
the eigenvalue $\lambda_Q=x+y\sqrt{D}$ of $M_Q$ is the
generator of the cyclic subgroup
\begin{displaymath}
\{z\in\mathcal{O}_D^\times\colon\mathcal{N}(z)=1\,\textit{and}\, z>1\}.
\end{displaymath}
Now consider that
the fundamental unit $\epsilon=x^\prime + y^\prime\sqrt{D}$
has negative norm $\mathcal{N}(\epsilon)=-1$.
Thus $\epsilon^2=\lambda_Q$
(since both, $\lambda_Q>1$ and $\epsilon>1$).
Hence $M\abcd{\epsilon}{0}{0}{\epsilon^{-1}}M^{-1}\in\textrm{SL}_2(K)$
is a square root of $M_Q$.
It is plain that
\begin{equation*}
M\mtwo{\epsilon}{0}{0}{\epsilon^{-1}}M^{-1}=
\frac{1}{\sqrt{D}}
\mtwo{Dy^\prime-Bx^\prime}{-2Cx^\prime}{2Ax^\prime}{Dy^\prime+Bx^\prime},
\end{equation*}
and the equality in the lemma follows.
Now using (once more) that $\epsilon>1$,
the latter part of the lemma follows.
\end{proof}

\begin{lem}
The matrix
\begin{displaymath}
S^{\frac{1}{2}}=\frac{1}{\sqrt{2}}\mtwo{1}{1}{-1}{1}
\end{displaymath}
is the square root of $S=\abcd{0}{1}{-1}{0}$;
it fixes $i$ and induces
the $z\mapsto iz$ map in the tangent space of $\hyper$ at $i$.
\end{lem}

\begin{proof}
The lemma follows after
a couple of straight forward calculations.
\end{proof}

\begin{defn}
Let $X$ be a topological space and
let $T_1$ and $T_2$ be positive real numbers.
Suppose we have
paths $f_1\colon[0,T_1]\longrightarrow X$
and $f_2\colon[0,T_2]\longrightarrow X$
such that $f_1(T_1)=f_2(0)$.
The \textit{concatenation} $f_1*f_2$ of $f_1$ and $f_2$ is the path
\begin{equation*}
f_1*f_2\colon[0,T_1+T_2]\longrightarrow X,
\end{equation*}
where
\begin{displaymath}
(f_1*f_2)(t)=\left\{
\begin{array}{ll}
f_1(t)&\textit{if $\,0\leq t<T_1$,}\\
f_2(t-T_1)&\textit{if $\,T_1\leq t\leq T_1+T_2$,}\\
\end{array}
\right.
\end{displaymath}
for each $t\in[0,T_1+T_2]$.
\end{defn}

\begin{thm}\label{thm:cont}
Assume that 
the base point $\tau_0$ of $g_{Q,\tau_0}$ is
not an elliptic point of order $e=2$.
Also,
if $g_{Q,\tau_0}$ contains an elliptic point of degree $e=2$,
further assume that the
fundamental unit $\epsilon$ of $\mathcal{O}_D$
has norm $\mathcal{N}(\epsilon)=-1$.
Then there is a
regular analytic path $\rho_{Q,\tau_0}$ on $X(1)$
such that,
up to a continuous reparametrisation of $g_{Q,\tau_0}$,
the restrictions of $\rho_{Q,\tau_0}$ and $g_{Q,\tau_0}$
to $I=[0,t_0)$ coincide
\begin{displaymath}
\rho_{Q,\tau_0}|I = g_{Q,\tau_0}|I,
\end{displaymath}
for $t_0\in\RR_{>0}$ sufficiently small.
Moreover,
up to a regular (analytic) reparametrisation
of $\rho_{Q,\tau_0}$,
the path $\rho_{Q,\tau_0}$
is uniquely determined by $g_{Q,\tau_0}$.
\end{thm}

\begin{proof}
Suppose $\gamma_Q$ contains
an elliptic point $\tau_e$ of order $e=3$.
By the commutative diagram
in the proof of Lemma~\ref{lem:reg}
we may see that for 
a small enough open neighbourhood $U$ of $\tau_e$
there is a continuous reparametrisation of
the restriction $g_{Q,\tau_0}|V$ of $g_{Q,\tau_0}$ to
the preimage $V=g_{Q,\tau_0}^{-1}U$,
such that $g_{Q,\tau_0}|V$ is both,
analytic and regular.
Now suppose $\gamma_Q$ contains
an elliptic point $\tau_e$ of order $e=2$.
Without loss of generality
we may assume $\tau_e=\sqrt{-1}\in\hyper$.
The proof splits in two steps as follows.

\noindent\textbf{ Step 1.} We claim that
the path $g_{Q,\tau_0}$ may be decomposed as
\begin{equation*} 
g_{Q,\tau_0}=
\{\tau_0,b\}_{\Gamma(1)}*\{b,\tau_0\}_{\Gamma(1)}*
\{\tau_0,b^\prime\}_{\Gamma(1)}*\{b^\prime,\tau_0\}_{\Gamma(1)},
\end{equation*}
for suitably chosen
elliptic points $b$ and $b^\prime$ in $\hyper$ of order $2$.
Indeed,
put $b=\tau_e$ and $b^\prime=M_Q^{-\frac{1}{2}}b$.
Now it suffices to prove that $b^\prime$ is an imaginary quadratic number
of discriminant $\Delta=-4$.
If we consider the positive definite form $P=[1,0,1]$,
our claim follows if we prove that
$P^\prime=M_Q^{-\frac{1}{2}}\circ P$
is integral,
primitive,
and has discriminant $\Delta=-4$.
From Definition~\ref{defn:circ} and Lemma~\ref{lem:sqr} we get
\begin{displaymath}
P^\prime=[1+2x^2-2Bxy, 8Axy, 1+2x^2+2Bxy],
\end{displaymath}
where $x+y\sqrt{D}\in\mathcal{O}_D^\times$
is a fundamental unit (of negative norm)
of the real quadratic order $\mathcal{O}_D$,
and $D$ is the discriminant of $Q$.
So $x$ and $y$ both lie in $\frac{1}{2}\ZZ$.
(By a slight abuse of notation
we have written $x$ and $y$ instead of $x^\prime$ and $y^\prime$.)
In particular we may see $P^\prime$ is integral.
Now suppose $P^\prime$ is not primitive.
So $M_Q^{-\frac{1}{2}}\circ P^\prime=\abcd{0}{-1}{1}{0}\circ P$
is not primitive,
but $P=[1,0,1]$ is obviously primitive.
We have reached a contradiction.
Therefore $P^\prime$ must be primitive.
Since it is plain that $P^\prime$ has discriminant $\Delta=-4$,
our claim follows.
A consequence of our claim is that
the path $g_{Q,\tau_0}$ not only fails to be regular at
the points $b$ and $b^\prime$;
these points are cusps
(in the sense of differential geometry)
of the path $g_{Q,\tau_0}$.

\noindent\textbf{Step 2.}
Now we ``remove'' the cusps of $g_{Q,\tau_0}$
by replacing in the above decomposition the components
\begin{displaymath}
\{b,\tau_0\}_{\Gamma(1)}*\{\tau_0,b^\prime\}_{\Gamma(1)}*\{b^\prime,\tau_0\}_{\Gamma(1)}
\end{displaymath}
by a certain concatenation of paths that meet at
an angle of $\frac{\pi}{2}$ at elliptic points of degree $e=2$.
More precisely,
we claim that
a suitable continuous reparametrisation of the concatenation
\begin{equation*}
p=
\{\tau_0,b\}_{\Gamma(1)}*
\{b,b^\prime\}_{\Gamma(1)}*
\{b^{\prime\prime},\tau_0\}_{\Gamma(1)},
\end{equation*}
is the desired regular path,
where
$b^\prime =M^{\frac{1}{2}}_{S^{\frac{a}{2}}\circ Q}b$
with $a=\textit{sign}(A)$,
and $b^{\prime\prime} =M_Q^{-\frac{1}{2}}b$.
Since we assumed that $\tau_0\not=b$,
the length $t_0=l(\tau_0,b)>0$.
So $p|I = g_{Q,\tau_0}|I$,
where $I=[0,t_0)$.
Again by the commutative diagram
in the proof of Lemma~\ref{lem:reg},
the fact that the ramification degree $e=2$ implies that
the path $p$ is both,
regular and analytic at the point $b$,
after a continuous reparametrisation of $p$.
Clearly the same is true for $p$ at $b^\prime$
(considering $b^\prime$ and $b^{\prime\prime}$
are the same point in $X(1)$).
Denote $\rho_{Q,\tau_0}$ any such  
continuous reparametrisation of $p$.
It is clear that $\rho_{Q,\tau_0}$ is unique
up to a regular analytic parametrisation
and the theorem follows.
\end{proof}

\begin{example}
Consider the indefinite form $Q=[1,4,-1]$ of discriminant $D=20$,
and let $\tau_0$ be the imaginary quadratic number
$\tau_0=\frac{-3+\sqrt{-19}}{2}\in\gamma_Q$.
Clearly $b=i$,
$b^\prime=i-4$,
and $b^{\prime\prime}=i+1$,
as depicted in Figure~\ref{cap:cycle0-20b}.
\begin{figure}
\begin{center}
\includegraphics{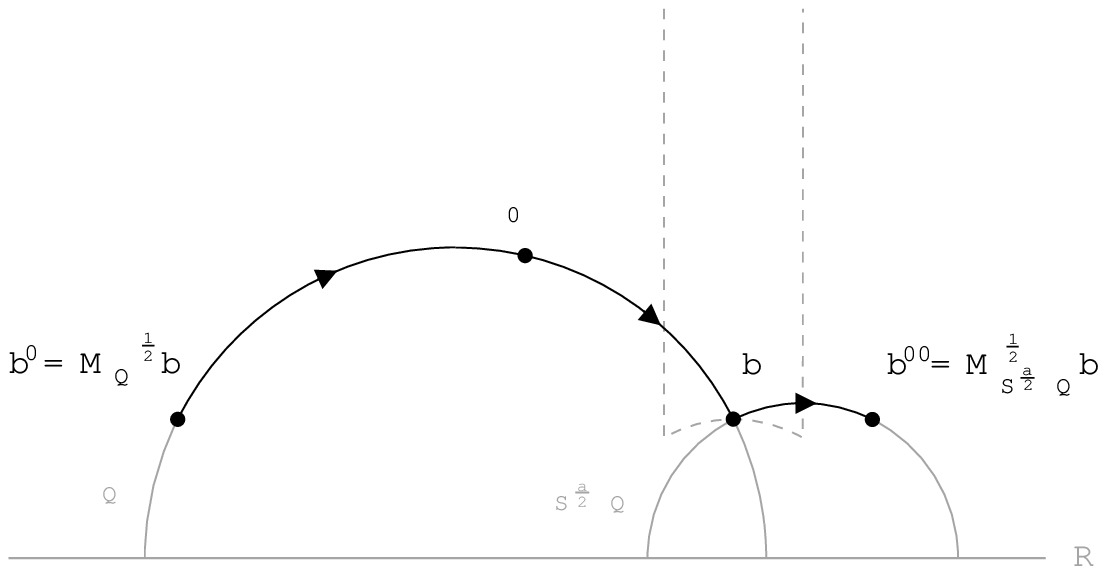}%
\end{center}
\caption{The path
$\{\tau_0,b\}*\{b,b^\prime\}*\{b^{\prime\prime},\tau_0\}$.}
\label{cap:cycle0-20b}
\end{figure}
Note the restriction $g_{Q,\tau_0}|I$
of $g_{Q,\tau_0}$ to  $I=[0,T)$ is regular,
where $T=l(\{\tau_0,b\})=1.216185939\dots>0$.\footnote{
It is a curious fact that
\begin{displaymath}
e^{l(\{\tau_0,b)\}}=\sqrt{\frac{218+96\sqrt{5}}{38}},
\end{displaymath}
which obviously lies in a quadratic extension of the field $\QQ(\sqrt{5})$.}
Note $g_{Q,\tau_0}$ fails to be regular at $b$.
Moreover,
the path $g_{Q,\tau_0}$ clearly has a cusp
(in the sense of differential geometry) at $b$.
Since the geodesics $\gamma_Q$ and $\gamma_{S^\frac{a}{2}\circ Q}$
meet at an angle of $\frac{\pi}{2}$,
and locally the function $j(\tau)$ at $b$
is the map $z\mapsto z^2$,
it is plain that
the path $\{b^\prime,b\}*\{b,b^{\prime\prime}\}$
is analytic and regular at $b$,
if $\{b^\prime,b\}*\{b,b^{\prime\prime}\}$
is suitably reparametrised.
\end{example}

\begin{rem}
As before suppose that
$\gamma_Q$ contains an elliptic point $\tau_e$ of order $e=2$.
Note that
the path $g_{Q,\tau_0}$ is generically $2:1$ onto its image,
whereas the path $\rho_{Q,\tau_0}$ is generically $1:1$ onto its image.
This explains the italics in the word ``essentially''
by the beginning of this subsection,
when describing how $g_{Q,\tau_0}$ traces out its image.
\end{rem}

%
\subsection{Geometric coding algorithm I}
Let $Q$ be an indefinite form and
let $\rho_{Q,\tau_0}$ be the
closed regular path we attached to $Q$ and $\tau_0\in\gamma_{Q}$
in the above subsection.
Here our main result is Algorithm~\ref{alg:cycle},
which produces
the (geometric) code of $\rho_{Q,\tau_0}$
with respect to
the fundamental region $\mathcal{F}$ of
the modular group $\Gamma(1)$.
So our algorithm may be regarded as a variant on
the geometric coding algorithm of Katok~\cite{katok:coding}.

We need to introduce some preliminaries and notations.
Let
\begin{displaymath}
S_Q=\left\{
\begin{array}{ll}
(S^{\frac{1}{2}})^a, &\textit{if $\,A+C=0$, where $a=\textrm{sign}(A)$}\\
S_{\pm 1},&\textit{if $\,2(A+C)=\pm B\not=0$}\\
S,&\textit{otherwise,}\\
\end{array}
\right.
\end{displaymath}
where
\begin{equation*}
S=\mtwo{0}{1}{-1}{0},\quad
S_{\pm 1}=\mtwo{\pm 1}{0}{1}{\pm 1},\quad\textrm{and}\quad
S^{\frac{1}{2}}=\frac{1}{\sqrt{2}}\mtwo{1}{1}{-1}{1}.
\end{equation*}
Note
\begin{equation*}
\sigma_{S^{\pm\frac{1}{2}}}=
\frac{1}{2}
\begin{pmatrix}
1& \pm 1& 1\\
\mp 2& 0& \pm 2\\
1& \mp 1& 1\\
\end{pmatrix},
\end{equation*}
where $M\mapsto\sigma_M$ is the homomorphism defined
by the beginning of this chapter.
So we have
\begin{displaymath}
S^{\pm\frac{1}{2}}\circ [A,B,-A]=
\left[\pm\frac{1}{2}B, \mp 2A, \mp\frac{1}{2}B\right]
\end{displaymath}
In particular,
$S^{\pm\frac{1}{2}}\circ Q$ may be non-primitive---and even non-integral.
Given any form $Q$ with coefficients in  $\QQ$,
we define $\nu_Q\in\QQ_{>0}$ to be such that $\nu_Q Q$
is integral and primitive. 

\begin{algorithm}\label{alg:cycle}
\KwData{a pair $(Q,\tau_0)$ with $\tau_0\gamma_Q$ as above}
\KwResult{the cycle $(Q)=(Q_0,Q_1,\dots, Q_l)$ of $Q$}
$Q \longleftarrow Q^\textit{red}$\;
$Q_0\longleftarrow Q$\;
$n \longleftarrow 0$\;
\While{$Q\not=Q_n$ or $n=0$}{
$Q_{n+1} \leftarrow (\nu_{Q_n} S_{Q_n}\circ Q_n)^\textit{nrm}$\;
$n\leftarrow n+1$\;
}
\caption{The cycle of $Q$}
\end{algorithm}

\begin{defn}
If $(Q)=(Q_0,Q_1,\dots Q_{l-1})$
is the cycle of $(Q,\tau_0)$
produced by Algorithm~\ref{alg:cycle},
the \textit{code} of the path $\rho_{Q,\tau_0}$ is
the sequence of matrices $(M_0,M_1,\dots, M_{l-1})$ such
that $Q_{n+1}=M_n\circ Q_n$,
obtained while performing Algorithm~\ref{alg:cycle}.
\end{defn}

\begin{example}
Consider the reduced form $Q=[1, 12, -1]$
of discriminant $D=2^2\cdot 37$.
The cycle and code of $Q$ are as follows.
\begin{displaymath}
\begin{array}{rlr}
n & Q_n & M_n\\
0 &[1, 12, -1] &\abcd{2}{0}{-1}{1}\\
1 &[3, 5, -1] &\abcd{5}{-1}{1}{0}\\
2 &[-1, 5, 3] &\abcd{-1}{-1}{1}{0}\\
3 &[3, 1, -3] &\abcd{13}{-11}{-1}{1}\\
\end{array}
\end{displaymath}
Note there is a ``jump'' in the discriminant
after we encounter a symmetric form.
This is a consequence of the fact that
the path $\gamma_Q$ contains
an elliptic point of order $e=2$.
\end{example}

Again let $Q$ be an indefinite form of discriminant $D$
(which we assume non-square).
It is seems likely that the condition $\mathcal{N}(\epsilon)=-1$
on the fundamental unit $\epsilon$ of $\mathcal{O}_D$
in our construction of the regular path $\rho_{Q,\tau_0}$
may be dropped. 
In fact,
Algorithm~\ref{alg:cycle}
makes perfect sense for any indefinite form $Q$
(with non-square, positive discriminant).
So one possible approach might be to
(1) define $\rho_{Q,\tau_0}$ as
the image in $X(1)$ of
a path constructed in some way out of
the code $(M_0,M_1,\cdots, M_{l-1})$ of $(Q,\tau_0)$
produced with the help of Algorithm~\ref{alg:cycle},
and then 
(2) prove that $\rho_{Q,\tau_0}^\prime$
(after a suitable continuous reparametrisation)
has the desired properties.
But we have not explored this approach any further
in this dissertation.
(In the applications we have in mind next chapter
it is possible to prove that $\mathcal{N}(\epsilon)=-1$.)

\begin{example}
Consider the reduced form $Q=[5, 9, -7]$ of discriminant $D=221$.
The cycle of the form $Q$ is as follows.
\begin{displaymath}
\begin{array}{rlr}
n & Q_n & M_n\\
0 &[5, 9, -7]& \abcd{1}{-1}{1}{0}\\
1 &[-7, 5, 7]& \frac{1}{2}\abcd{7}{5}{1}{1}\\
2 &[-5, 32, -7]& \abcd{4}{-1}{1}{0}\\
3 &[-7, 24, 11]& \abcd{-2}{-1}{1}{0}\\
4 &[11, 20, -11]& \abcd{3}{-1}{-1}{1}\\
\end{array}
\end{displaymath}
It is well-known that the fundamental unit of $\mathcal{O}_D$
has positive norm.
\end{example}

%
%
%
\section{Higher levels}\label{sec:leveln}

%
\subsection{Reduction algorithms}
Let $\Gamma$ be a group acting on a Riemann surface $X$.
A \textit{fundamental domain} for $\Gamma$ is
a connected open subset $D$ of $X$ such that the following two
properties hold.
\begin{description}
\item[FD 1] No two points on $D$ lie in the same orbit of $\Gamma$.
\item[FD 2] The closure $D^\textit{clo}$ of $D$ contains
at least one element from each orbit.
\end{description}
The \textit{fundamental region} associated to a fundamental domain $D$ is
the closure $D^\textit{clo}$ of $D$.
A \textit{fundamental set} may be defined as a set
$D\subset S\subset D^\textit{clo}$ such that the canonical map from 
$S$ into $\Gamma\backslash X$ is a bijection.

\begin{defn}
Suppose $S$ is a fundamental set for the action of $\Gamma^0(N)$
on $\hyper$.
We say a positive definite form $Q$ is 
\textit{$S$-reduced} if $\pi_Q\in \mathcal{F}_N$.
\end{defn}

We shall now describe a fundamental set for $\Gamma^0(N)$.
As before assume that $N$ is an odd prime.
Define
\begin{equation*}
\mathcal{R}(N)=\left\{-\frac{N-1}{2},\dots,-2,-1\right\}\bigcup
\left\{1,2,\dots,\frac{N-1}{2}\right\}
\end{equation*}
and let $T=\abcd{1}{1}{0}{1}$ and $S=\abcd{0}{-1}{1}{0}$
be the standard generators of $\textrm{SL}_2(\ZZ)$.
It is well-known that the union
\begin{equation*}
\mathcal{F}_N =
\{S\mathcal{F}^\textit{clo}\} \cup 
\left\{T^r\mathcal{F}^\textit{clo}\colon r\in\mathcal{R}(N) \right\},
\end{equation*}
is a fundamental region for $\Gamma^0(N)$ with an even number of sides
which are identified by the generators of $\Gamma^0(N)$
\begin{equation*}
T^N=\mtwo{1}{N}{0}{1},\quad\textrm{and}\quad
M_s=\mtwo{s^\prime}{-ss^\prime -1}{1}{-s},\quad
\in\Gamma^0(N),
\end{equation*}
with $s,s^\prime\in\mathcal{R}(N)$ such that $r^\prime r\equiv -1\pmod{N}$.
(See e.g. Apostol~\cite{apostol:modular}, p. 76).

\begin{figure}
\begin{center}
\includegraphics{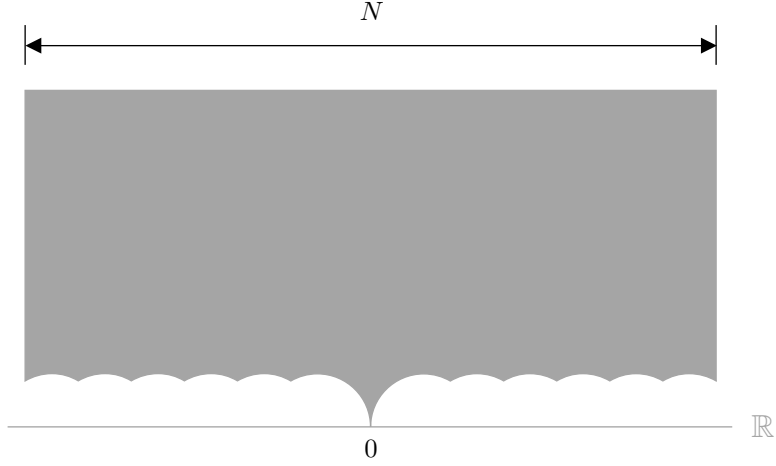}
\end{center}
\caption{Fundamental domain for $\Gamma^0(N)$, where $N=13$.}
\label{cap:leftfundom37}
\end{figure}

Now we briefly review a construction of a fundamental set
$\mathcal{F}_N^\prime$ such that
$\mathcal{F}_N=(\mathcal{F}_N^\prime)^\textit{clos}$.
First note that for each point
$\tau\in\partial(\mathcal{F}_N)=\mathcal{F}_N - \mathcal{F}_N^\textit{int}$
the orbit $\Gamma^0(N)\tau$ consists of at most $3$ points, say,
$\tau_1, \tau_2$, and $\tau_3$.
Clearly $\Re(\tau_i)\not=\Re(\tau_j)$ if $i\not=j$, for each $i,j = 1,2,3$.
So it makes sense to define the 
\textit{left most} point $l(\tau)$ of the orbit $\Gamma^0(N)\tau$
as the point $\tau_i$ with the least real part $\Re(l(\tau))\leq\Re(\tau_i)$,
$i= 1,2,3$.
Then we may define $\mathcal{F}_N^\prime$ as
\begin{equation*}
\mathcal{F}_N^\prime = \mathcal{F}_N^\textit{int} \cup
\{l(\tau):\tau\in\partial(\mathcal{F}_N)\}
\end{equation*}
which is by construction a fundamental set for $\Gamma^0(N)$ such that
$\mathcal{F}_N=(\mathcal{F}_N^\prime)^\textit{clos}$.

So given any positive definite form 
$Q$ there is an element $R_{Q,N}\in\Gamma^0(N)$ such that
$R_{Q,N}\circ Q$ is $\mathcal{F}_N$-reduced.
Such element may be computed as follows.
Using Algorithm~\ref{alg:gauss}
we may find a matrix  $M_1\in\textrm{SL}_2(\ZZ)$ such that
$Q$ is $\mathcal{F}$-reduced.
Now by trial and error we may pick a matrix $M_2$ among
the elements of the (finite) set of coset representatives
$\{S\}\cup \{T^n\colon n\in\mathcal{R}(N)\}$
of $\Gamma^0(N)$ in SL$_2(\ZZ)$ such that $M=M_2M_1\in\Gamma^0(N)$.
Note that the imaginary quadratic 
$\tau=\pi_{Q^\prime}$ attached to $Q^\prime=M\circ Q$ clearly lies in 
the fundamental region $\mathcal{F}_N$.
If $\tau\in\mathcal{F}_N^\textit{int}$ then 
$R_{Q,N}=M$ has required property.
So we suppose $\tau\in\partial(\mathcal{F}_N)$.
If $\tau\not\in\mathcal{F}_N^\prime$ then
for suitable choice among the generators $G$ of $\Gamma^0(N)$
introduced above we have
$G\tau$ lies in the fundamental set $\mathcal{F}_N^\prime$.
So $R_{Q,N}=GM$ has the required property and we are done.

Now we turn to the $D>0$ case.
Let $Q=[A,B,C]$ be an indefinite binary quadratic form
of discriminant $D = B^2 - 4AC>0$.
As before we assume $D$ is non-square,
to avoid forms $Q$ such that 
the geodesic $\gamma_{Q}$
contains a side of the fundamental region $\mathcal{F}_N$, e.g.
$Q=[1, 2, 0], [1, 4, 3], [1, 6, 8],\dots$.
Here we discuss an extension of the definition of
normalisation and reduction of $Q$ for the group $\Gamma^0(N)$.

\begin{lem}\label{lem:nrm}
Pick a base point $\tau_0$ of $\gamma_Q$,
e.g. $\tau_0=\widehat{\gamma_Q}$.
The form $Q$ is $\Gamma^0(N)$-equivalent to a form,
we denote $Q^\textit{mt}$, such that
$\gamma_{Q^\textit{mt}}$ meets $\mathcal{F}_N^\textit{int}$
\begin{equation*}
\gamma_{Q^\textit{mt}}\cap\mathcal{F}_N^\textit{int}\not=\emptyset.
\end{equation*}
\end{lem}

\begin{proof}
Now let $R_{\tau_0,N}$ denote the matrix in $\Gamma^0(N)$ 
such that $\tau_0^\prime = R_{\tau_0,N}\circ\tau_0$
is $\mathcal{F}_N$-reduced. 
Clearly the geodesic $\gamma_{Q_1}$
attached to $Q_1 = R_{\tau_0,N}\circ Q$
is such that
\begin{displaymath}
\gamma_{Q_1}\cap\mathcal{F}_N^\textit{int}\not=\emptyset,
\end{displaymath}
except possibly when $\tau_0^\prime = \rho_{\pm 1}\pm r$,
with $r\in\{1,2,\dots,\frac{N-1}{2}\}$,
i.e. when the point $\tau_0^\prime$ is one of the non-cusp vertices
of the fundamental region $\mathcal{F}_N$.
So let us suppose $\rho_{\pm 1}\pm r \in \gamma_{Q_1}$ and 
$\gamma_{Q_1}\cap \mathcal{F}_N^\textit{int}=\emptyset$.\footnote{
This case occurs very seldom in practice.
For example,
fix $b\in\RR_{>0}$ moderately large,
and also fix a positive (non-square) discriminant $D$.
Experimental evidence suggests that
the probability that a form $Q=[A,B,C]$ picked randomly in the set
\begin{displaymath}
\{Q=[A,B,C]\in\ZZ^3\,\colon\,
|A|<b,|B|<b,|C|<b,B^2-4AC=D,\,\textit{and $Q$ primitive}\}
\end{displaymath}
falls in this case is small.}
Put $Q_2=T^{-s}\circ Q_1$, where $s = \pm r$ and note that
$\gamma_{Q_2}$ contains the point $\rho_{\pm 1}$.
Now by Lemma~\ref{lem:rhoint} we know that
\begin{equation*}
Q_3=S\circ Q_2=[C_2,-B_2,A_2]
\end{equation*}
is such that $\gamma_{Q_3}$
has non-trivial intersection with the open set $\mathcal{F}^\textit{int}$.
Therefore $Q_4 = T^{s^\prime}\circ Q_3$ is such that
$\gamma_{Q_4}$ has non-trivial intersection with
$T^{s^\prime}\mathcal{F}^\textit{int}\subset\mathcal{F}_N^\textit{int}$,
for each $s^\prime\in\mathcal{R}(N)$.
In particular if $s s^\prime \equiv -1 \pmod{N}$ then the form
\begin{equation*}
Q_4 = (T^{s^\prime}S T^{-s})\circ Q_1 =
(T^{s^\prime}S T^{-s} R_{\widehat{Q},N})\circ Q
\end{equation*}
is clearly $\Gamma^0(N)$-equivalent to $Q$,
and it is plain that
$\gamma_{Q_4}\cap\mathcal{F}_N^\textit{int}\not=\emptyset$.
\end{proof}

Given a geodesic segment $\gamma=\{\alpha,\beta\}$
let $\gamma^\textit{int}$ be the open geodesic segment
obtained by removing the endpoints $\alpha$ and $\beta$ of $\gamma$.
So let $\sigma_{\pm 1,\infty}\subset\partial(\mathcal{F}_N)$
be the open geodesic segment defined by
\begin{equation*}
\begin{split}
\sigma_{-1,\infty}&=\left\{\infty,\rho_{-1}-\frac{N - 1}{2}\right\}^\textit{int}\\
\sigma_{+1,\infty}&=\left\{\rho_{+1}+\frac{N - 1}{2},\infty\right\}^\textit{int}
\end{split}
\end{equation*}
In a similar way let us define the open geodesic segment
$\sigma_{\pm 1,0}\subset\partial(\mathcal{F}_N)$ by
\begin{equation*}
\begin{split}
\sigma_{-1,0}&=\left\{\rho_{-1}, 0\right\}^\textit{int}\\
\sigma_{+1,0}&=\left\{0,\rho_{+1}\right\}^\textit{int}
\end{split}
\end{equation*}
Now for the rest of this paragraph
let us restrict the discussion only to
forms $Q$ such that 
$\gamma_Q\cap\mathcal{F}_N^\textit{int}\not=\emptyset$.
We say $Q$ is \textit{close to the cusp $\infty$} if
$I(\gamma_Q, \sigma_{\textit{sign}(A),\infty}) = 1$.
This is equivalent to the condition
\begin{equation*}
3|A| + \textit{sign}(A)(4C + N^2A) + 2NB < 0.
\end{equation*}
The \textit{$N$-normalisation} $Q^{N\textit{-nrm}}$
of a form $Q$ close to the cusp $\infty$
is the $\Gamma^0(N)$-equivalent form
$Q^{N\textit{-nrm}}=T^{\left[\frac{\delta(Q)}{N}\right]N}\circ Q$,
where $[\cdot]$ is the nearest integer function.
We say $Q$ is \textit{close to the cusp $0$} if
$I(\gamma_Q, \sigma_{\textit{sign}(A),0}) = 1$.
This means that
\begin{equation*}
3|C| + \textit{sign}(C)(4A + C) - 2B < 0.
\end{equation*}
The \textit{$N$-normalisation} $Q^{N\textit{-nrm}}$
of a form $Q$ close to the cusp $0$ is the 
$\Gamma^0(N)$-equivalent form
$Q^{N\textit{-nrm}}=(S T^{\delta(S\circ Q)} S)\circ Q$,
where $\delta$ is as in Definition~\ref{defn:nrm}.
For completeness,
if $Q$ is neither close to the cusp $\infty$
nor close to the cusp $0$ we define the
\textit{$N$-normalisation} $Q^{N\textit{-nrm}}$ of $Q$ to be $Q$ itself.
Finally,
the $N$-\textit{reduction} $Q^{N\textit{-red}}$ of an arbitrarily given
indefinite form $Q$ (with non-square discriminant) is the
normalisation
$Q^{N\textit{-red}}=(Q^\textit{mt})^{N\textit{-nrm}}$
of $Q^\textit{mt}$,
where $Q^\textit{mt}$ is the form defined in Lemma~\ref{lem:nrm}.
If a form $Q_1$ is the $N$-reduction of a form $Q_0$ then we say
$Q_1$ is $N$\textit{-reduced}.

\begin{figure}
\begin{center}
\includegraphics{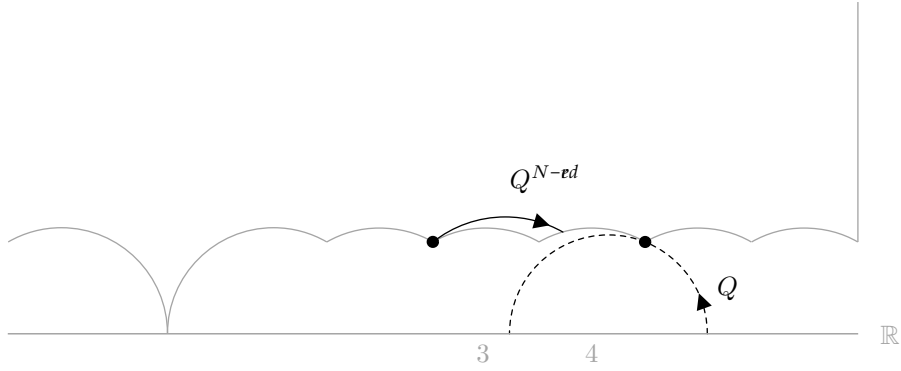} %
\end{center}
\caption{The $N$-reduction of $Q=[-13,108,-213]$, with $N=13$.}
\label{cap:red013}
\end{figure}

\begin{example}
Let $N = 13$ and consider the primitive
indefinite form
\begin{equation*}
Q = [-13, 108, -213]
\end{equation*}
of discriminant $D=588$.
Note $\gamma_Q\cap\mathcal{F}_N^\textit{int}=\emptyset$.
Also $\rho_{+1}+r\in\gamma_Q$ for $r=4$.
So by the proof of Lemma~\ref{lem:nrm} we may see that
$Q^\textit{mt} = (T^{3}ST^{-4})\circ Q$
is a $\Gamma^0(N)$-equivalent form such that
$\gamma_{Q^\textit{mt}}$ meets the open set $\mathcal{F}_N^\textit{int}$.
Now since $Q^\textit{mt}$ is neither
close to the cusp $\infty$ nor close to the cusp $0$,
then by definition $Q^{N\textit{-red}}=Q^\textit{mt}$,
as depicted in Figure~\ref{cap:red013}.
\end{example}

%
\subsection{Geometric coding algorithm II}\label{sub:codii}
Note that $\Gamma^0(N)=S^{-1}\Gamma_0(N)S$,
where $S=\abcd{0}{-1}{1}{0}$.
By a slight abuse of notation\footnote{
It is more customary to
write $X^0(N)=\Gamma^0(N)\backslash\hyper^*$.
But we want to avoid writing
non-standard expressions such as $X_+^0(N)$ (or $X^{0+}(N)$)
when considering the Atkin-Lehner quotient defined by
the Fricke involution $w_N$.
}
from now on we write
\begin{displaymath}
X_0(N)=\Gamma^0(N)\backslash\hyper^*.
\end{displaymath}
As before let $Q=[A, B, C]$ be primitive indefinite form
with non-square discriminant $D=B^2 - 4AC$.
Now define $M_{N,Q}=M_Q^n$,
where 
\begin{equation*}
M_Q = \mtwo{u-Bv}{-2Cv}{2Av}{u+Bv}\in\Gamma(1),
\end{equation*}
and $(u,v)\in\ZZ$ is a
fundamental solution of the ordinary Pell equation
\begin{equation*}
X^2-DX^2 = 1,
\end{equation*}
and $n$ is the smallest positive integer such that
$M_Q^n\in\Gamma_0(N)$.
Now pick a base point $\tau_0$ of $\gamma_Q$ and
let $g_{N; Q,\tau_0}$ denote
the image of the path $\{\tau_0, M_{N, Q}\tau_0\}$
in the modular curve $X_0(N)$.
The construction
of the regular path
in the proof of Theorem~\ref{thm:cont}
carries over to $X_0(N)$ as follows.
It is well-known that
the elliptic points $\tau_e\in\hyper$
for the action of $\Gamma^0(N)$ in the upper half plane $\hyper$ are
precisely the ones attached to
primitive forms $P=[A,B,C]$ with $N\mid C$ and
discriminant $\Delta=-3$,
if the degree of the point is $e=3$ and
discriminant $\Delta=-4$,
if the degree of the point is $e=2$.
Assume there is an
elliptic point $\tau_e$ of degree $e=2$,
e.g. $\tau_e = r + \sqrt{-1}\in\hyper$
with $r\in\ZZ$ such that
\begin{displaymath}
r^2\equiv -1\pmod{N},
\end{displaymath}
and that $\gamma_Q$ contains such elliptic point $\tau_e$.
The matrix $S^\frac{1}{2}$ in the said construction
is to be replaced by $T^r S^\frac{1}{2} T^{-r}$,
where and $T=\abcd{1}{1}{0}{1}$.
Again let $n$ be the smallest positive integer
such that $M_Q^n$ lies in $\Gamma^0(N)$.
Further assume that
the fundamental unit $\epsilon$ of
the real quadratic order $\mathcal{O}_D$
has norm $\mathcal{N}(\epsilon)=-1$. 
The matrix $M_Q^\frac{1}{2}$ is to be
replaced by $M_{N,Q}^\frac{1}{2}=(M_Q^\frac{1}{2})^n$,
so that $(M_{N,Q}^\frac{1}{2})^2=M_{N,Q}$.
Let $\rho_{N, Q,\tau_0}$ be
the regular analytic path thus obtained
and call it the \textit{regular path of level $N$} attached to $Q$.
This path is unique up to a regular analytic reparametrisation.

Now we need to extend Definition~\ref{defn:nrm} to include
some indefinite forms $Q$ that arise
in Algorithm~\ref{alg:ncycle}.
As before put $\sigma=\{\rho_{-1},\rho_{+1}\}$.
Let $\delta(Q)$ be just as in Definition~\ref{defn:nrm},
if $D > 4A^2$ and let $\delta(Q)$ be the integer such that
the intersection number $I(\gamma_Q + \delta(Q),\sigma) = 1$
(whenever it exists), if $3A^2D < 4A^2$.
Assume that $Q$ is $N$-reduced
and that $Q$ is neither close to cusp $\infty$ nor close to cusp $0$.
Now put $a = \textit{sign}(A)$, $d = \delta(Q)$,
and define the matrix
\begin{equation}\label{eqn:sqn}
S_{N,Q}=\left\{
\begin{array}{ll}
T^a S T^{-a}, &
\textit{if $\,d = 0$, i.e. $\rho_a\in\gamma_Q$}
\quad\hfill\textrm{(1)}\\
T^{-a\frac{N + 1}{2}} S T^{-2a}, &
\textit{if $\,\rho_a \in\gamma_{Q} - a$}
\quad\hfill\textrm{(2)}\\
T^{(d-a)^*} S T^{d-a}, &
\textit{if $\,\rho_a \in\gamma_{Q} + d$ and $|d|\not = 1$}
\quad\hfill\textrm{(3)}\\
T^{d^*} S^{\frac{a}{2}} T^d, &
\textit{if $\,i_N\in\gamma_{Q}$}
\quad\hfill\textrm{(4)}\\
T^{d^*} S T^d, &
\textit{otherwise,}
\quad\hfill\textrm{(5)}\\
\end{array}
\right.
\end{equation}
where $i_N = i + k$ and $k\in C(N)$ is such that
$k^2\equiv -1\pmod{N}$,
whenever $N\equiv 1\pmod{4}$,
and $s^*$ is the element of $C(N)$ such that
$s s^* \equiv 1\pmod{N}$,
given $s\in\ZZ$ with $s\not\equiv 0\pmod{N}$. 

\begin{algorithm}\label{alg:ncycle}
\KwData{a primitive indefinite binary quadratic form $Q$}
\KwResult{the $N$-cycle $(Q)_N = (Q_0,Q_1,\dots, Q_{l-1})$ of $Q$}
$Q \longleftarrow Q^{N\textit{-red}}$\;
$Q_0\longleftarrow Q$\;
$n \longleftarrow 0$\;
\While{$Q\not=Q_n$ or $n=0$}{
\eIf{$Q_n$ is either near cusp $\infty$ or cusp $0$}{$Q_{n+1}\longleftarrow
Q_n^{N\textit{-nrm}}$}{$Q_{n+1}\longleftarrow\nu_{Q_n}S_{N,Q_n}\circ Q_n$}
$n \longleftarrow n + 1$\;
}
\caption{The $N$-cycle of $Q$}
\end{algorithm}

Now we shall consider some examples
of $N$-cycles $(Q)_N=(Q_0,\dots,Q_{l-1})$
that may help to clarify the above ideas.
In the tables below the 
last column encodes the steps followed while executing
the algorithm,
using the the convention that
a number from $1$ to $5$ indicates one of the five cases
of the definition of matrix $S_{N,Q}$,
while $0$ (resp. $\infty$) indicates the normalisation step
associated to a
form close to the cusp $0$ (resp. the cusp $\infty$).

\begin{example}
Suppose $N=13$.
Consider the indefinite form $Q = [-13, 108, -213]$
of discriminant $D=7^2\cdot 12$.
The $N$-reduction of $Q^{N\textit{-red}}=[11, -70, 98]$ of $Q$
has the following $N$-cycle.
\begin{displaymath}
\begin{array}{rllr}
n & Q_n & M_n & \textrm{Case}\\
0 &[11, -70, 98]& \abcd{3}{-13}{1}{-4}&\textrm{5}\\
1 &[-6, 18, 11]& \abcd{1}{0}{1}{1}& \textrm{0}\\
2 &[-13, -4, 11]& \abcd{1}{0}{1}{1}&\textrm{5}\\
3 &[2, -26, 11]& \abcd{1}{-13}{0}{1}&\textrm{$\infty$}\\
4 &[2, 26, 11]& \abcd{1}{0}{1}{1}&\textrm{5}\\
5 &[-13, 4, 11]& \abcd{1}{0}{1}{1}&\textrm{0}\\
6 &[-6, -18, 11]& \abcd{-4}{-13}{1}{3}&\textrm{5}\\
7 &[11, 70, 98]& \abcd{-6}{-13}{1}{2}&\textrm{3}\\
8 &[2, -2, -73]& \abcd{-2}{13}{1}{-7}&\textrm{3}\\
9 &[11, 18, -6]& \abcd{1}{0}{-3}{1}&\textrm{0}\\
10 &[11, -18, -6]& \abcd{-7}{13}{1}{-2}&\textrm{2}\\
11 &[2, 2, -73]& \abcd{2}{-13}{1}{-6}&\textrm{3}\\
\end{array}
\end{displaymath}
Figure~\ref{cap:cycle013a} depicts the intersection of
the geodesics $\gamma_{Q_n}$,
with the fundamental domain $\mathcal{F}_N$,
for $n=0,1,2$ and $3$.
Similarly,
Figure~\ref{cap:cycle013b} depicts
the intersection of the geodesics $\gamma_{Q_n}$,
with the fundamental domain $\mathcal{F}_N$,
for $n=9,10$ and $11$.
Note that instance (4) of $S_{N,Q}$ in Equation~\ref{eqn:sqn}
does not take place.
This means that the geodesic $\gamma_Q$ does not contain
an elliptic point $\tau_e$ of degree $e=2$.
\begin{figure}
\begin{center}
\includegraphics{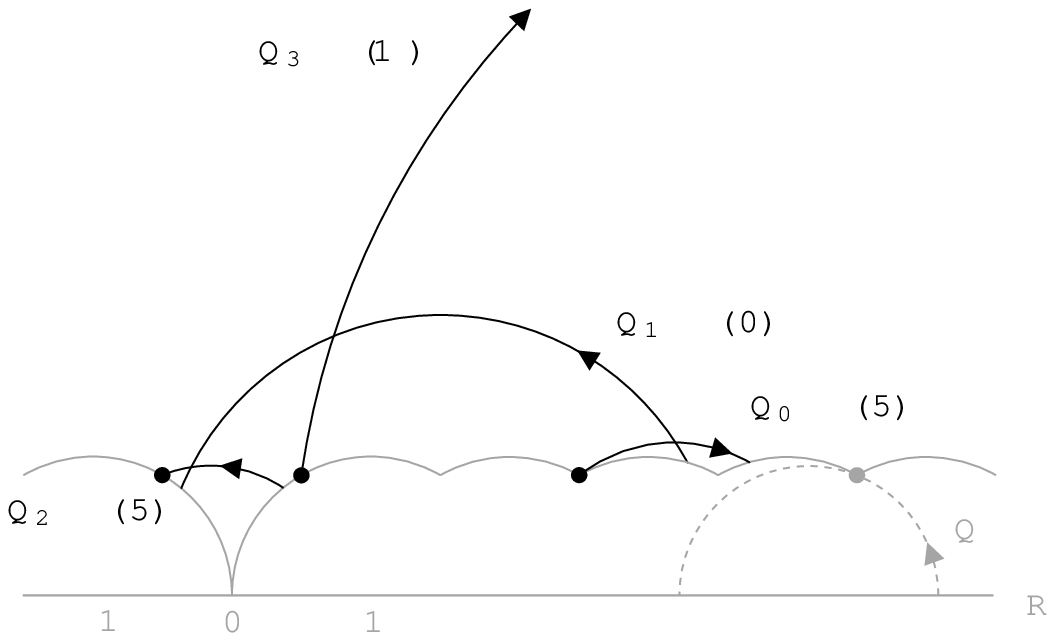} %
\end{center}
\caption{Beginning of the $N$-cycle of $Q=[-13,108,-213]$, with $N=13$.}
\label{cap:cycle013a}
\end{figure}
\begin{figure}

\begin{center}
\includegraphics{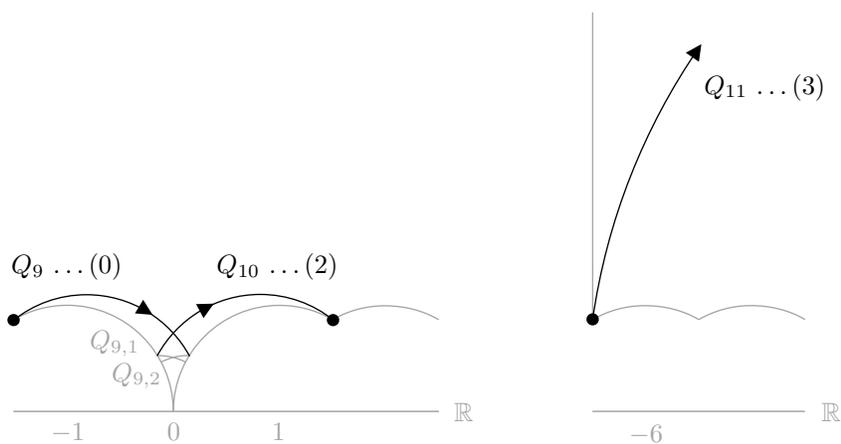} %
\end{center}
\caption{Near the end of the $N$-cycle of $Q=[-13,108,-213]$, with $N=13$.}
\label{cap:cycle013b}
\end{figure}
\end{example}

\begin{example}\label{ex:5-13}
Now suppose $N= 5$ and consider
the $N$-reduced form $Q=[1, -1, -3]$ of discriminant $D=13$.
The $N$-cycle of $Q$ is as follows.
\begin{displaymath}
\begin{array}{rllr}
n & Q_n & M_n & \textrm{Case}\\
0 &[1, -1, -3]& \abcd{-1}{5}{-1}{3}&\textrm{4}\\    
1 &[3, -16, 17]& \abcd{1}{-5}{0}{1}&\textrm{$\infty$}\\  
2 &[3, 14, 12]& \abcd{1}{0}{1}{1}&\textrm{3}\\  
3 &[1, -10, 12]& \abcd{1}{-10}{0}{1}&\textrm{$\infty$}\\    
4 &[1, 10, 12]& \abcd{1}{0}{1}{1}&\textrm{3}\\  
5 &[3, -14, 12]& \abcd{1}{-5}{0}{1}&\textrm{$\infty$}\\   
6 &[3, 16, 17]& \abcd{3}{5}{-1}{-1}&\textrm{4}\\    
7 &[1, 1, -3]& \abcd{-1}{0}{1}{-1}&\textrm{5}\\    
8 &[-1, -5, -3]& \abcd{1}{5}{0}{1}&\textrm{$\infty$}\\    
9 &[-1, 5, -3]& \abcd{-1}{0}{1}{-1}&\textrm{5}\\    
\end{array}
\end{displaymath}
Figure~\ref{cap:cycle013-005} depicts the intersections of
the geodesics $\gamma_{Q_n}$
with the fundamental domain $\mathcal{F}_N$,
for $n=1,2$, and $3$.
Note that $Q_n$ falls in case (4),
for $n=0$ and $6$.
This means that $\gamma_Q$
contains an elliptic point $\tau_e$ of order $e=2$.
\begin{figure}
\begin{center}
\includegraphics{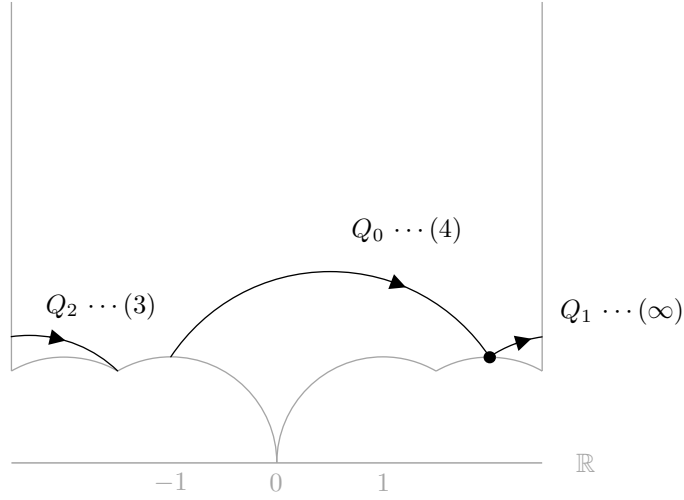} %
\end{center}
\caption{Beginning of the $N$-cycle of $Q=[1, -1, -3]$, with $N=5$.}
\label{cap:cycle013-005}
\end{figure}
\end{example}

\begin{rem}
It is not hard to see that
for each prime $N>1$,
the $N$-cycle of the the indefinite form $Q=[2, 2N, -1]$
is given by
\begin{displaymath}
\begin{array}{rlr}
n & Q_n & M_n\\
0 &[2, 2N, -1] & \abcd{1}{0}{-2N}{1}\\
1 &[2,-2N,-1]  & \abcd{1}{-N}{0}{1}\\
\end{array}
\end{displaymath}
Note that
the matrices $M_0=\abcd{1}{0}{-2N}{1}$ and $M_0=\abcd{1}{-N}{0}{1}$
are parabolic.
So the code $(M_0,M_1,\dots)$
of some indefinite forms $Q$
may consist of parabolic elements only. 
\end{rem}

\begin{defn}
The $N$-\textit{code} of $Q$ is the sequence of matrices
$(M_0,M_1,\dots)$ defined by $Q_{n+1}=M_n\circ Q_n$,
where the matrices $M_n$ come from the computation of
the $N$-cycle $(Q)_N=(Q_0,Q_1,\dots Q_l)$ of $Q$ using
Algorithm~\ref{alg:ncycle}.
\end{defn}


%
%
%
\bibliographystyle{amsplain}
\bibliography{biblio}
\printindex
\end{document}